\documentclass[12pt]{amsart} 
\usepackage{amssymb}
\usepackage{euscript}

\makeatletter \@mparswitchfalse \makeatother

\newtheorem{theorem}{Theorem}
\newtheorem{proposition}[theorem]{Proposition}

\newtheorem{corollary}[theorem]{Corollary}

\begin{document}

% bring back \eqalign from plain TeX
\def\eqalign#1{\null\,\vcenter{\openup\jot
  \ialign{\strut\hfil$\displaystyle{##}$&$\displaystyle{{}##}$\hfil
      \crcr#1\crcr}}\,}

\let\eps\varepsilon
\def\R{\mathbb R}
\let\xi\zeta
\let\sharp\#
\let\cal\EuScript
\let\\\cr
\let\phi\varphi
\let\union\bigcup
\let\inter\bigcap
\let\kappa=k
\def\supp{\operatorname{supp}}
\def\loc{\operatorname{loc}}
\def\fin{\operatorname{fin}}
\let\blip\rho
\let\pee=P
\let\emptyset\varnothing
\let\curly=r
\let\jay=j
\let\capitalu=U
\let\lceil\Gamma

\title[Algebras associated with Blaschke products
of type $G$]{Algebras associated with Blaschke
products of type $\boldsymbol G$}
\author{Carroll Guillory}
\address{\hskip-\parindent Carroll Guillory\\
Department of Mathematics\\
University of Southwestern Loui\-si\-ana\\
Lafayette, LA 70504}
\email{cjg2476@@usl.edu}

\author{Kin Y. Li}
\address{\hskip-\parindent 
Kin Y. Li\\
Department of Mathematics\\ 
Hong Kong University of Science and Technology\\ 
Clear Water Bay, Kowloon, Hong Kong}
\email{makyli@@uxmail.ust.hk}

\begin{abstract}
Let $\Omega$ and $\Omega_{\fin}$ be the sets of all interpolating
Blaschke products of type $G$ and of finite type $G$, respectively.  
Let $E$ and $E_{\fin}$ be the Douglas algebras 
generated by $H^\infty$ together with the complex conjugates of elements of 
$\Omega$ and $\Omega_{\fin}$, respectively.
We show that the set of all invertible inner
functions in $E$  is the set of all finite products 
of elements of $\Omega$ , which is also the closure of $\Omega$  among the Blaschke products. 
Consequently, finite convex combinations of
finite products of elements of $\Omega$   are
dense in the closed unit ball of the subalgebra of $H^\infty$ generated by
$\Omega$.  The same results hold when we replace $\Omega$ by
$\Omega_{\fin}$ and $E$ by $E_{\fin}$.
\end{abstract}

\maketitle

\section{Introduction} Let $D$ be the open unit disk and $T$ be the unit circle on the complex plane. Let $H^\infty$ be the Banach algebra of bounded
analytic functions on the open unit disk $D$. Via radial limits we
can consider $H^\infty$ as a closed subalgebra of $L^\infty$,
 where $L^\infty$ is the family of all essential bounded measurable
functions  on $T$. 
Any function $h$ in $H^\infty$ with $|h|=1$ {\it a.e.} on $T$ is called an {\em inner} function.
Let $\{z_n\}$ be a sequence in $D$ with $\sum_n(1-|z_n|) <
\infty$. Then the function 
$$
b(z)=\prod_n  \frac{\bar{z}_n}{|z_n|}
\frac{z_n-z}{1-\bar{z}_n z}\quad\hbox{for $z\in D$,}
$$
is called a {\em Blaschke product} with roots $\{z_n\}$. Let
$$
\delta (b) = \inf_{\kappa} \prod_{n\neq \kappa}\left | \frac{z_\kappa
- z_n}{1- \bar{z}_n z_\kappa }\right |.
$$
If $\delta(b) > 0$, then $b$ and $\{z_n\}$ are called {\em interpolating}.
By \cite{Carl1}, if $\delta(b) > 0$, then for every bounded sequence
$\{a_n\}$, there exists $f$ in $H^\infty$ such that $f(z_n)=a_n$ for
every $n$. If
$$
\lim_{k\to \infty}\prod_{n\neq \kappa} \left | \frac{z_\kappa - z_n}
{1-\bar{z}_n z_\kappa} \right | = 1,
$$
then $b$ and $\{z_n\}$ are called {\em thin} or {\em sparse}.

We denote by
$M(H^\infty)$ the maximal ideal space of $H^\infty$. 
A closed subalgebra $B$ between $H^\infty$ and $L^\infty$ is called a
{\em Douglas algebra.} We denote by $M(B)$ the maximal ideal space of the
Douglas algebra $B$. For an interpolating Blaschke product $b$, we
denote by $H^\infty[\bar{b}]$ the Douglas algebra generated by
$H^\infty$ and the complex conjugate of $b$. For a function $f$ in
$H^\infty$, let
$$
Z(f)=\{x\in M(H^\infty): f(x) = 0\}
$$
and for $0<c\leq 1$,
$$
\{|f| < c\} = \{x \in M(H^\infty): |f(x)| < c\}.
$$
 For a point $x$ in $M(H^\infty)$, there is
a representing measure $\mu_x$ on $M(L^\infty)$, that is,
$$
f(x)=\int_{M(L^\infty)} f \,d \mu_x
$$
for every $f\in H^\infty$. We denote by  $\supp\mu_x$ the support set
for the representing measure $\mu_x$. 

By the Corona
Theorem, $D$ can be considered as a dense subset of $M(H^\infty)$. For
points $x,y$ in $M(H^\infty)$, let
$$
\blip (x,y) = \sup\{|f(y)|: f \in H^\infty , \| f\|_\infty \leq 1, f(x)
= 0\}
$$
and put
$$
\pee (x) = \{m\in M(H^\infty): \blip(m,x) < 1\}.
$$
The set $\pee(x)$ is called the {\em Gleason part containing} $x$. For $z, w  \in D$, we have
$$\blip(z,w ) = \left|\frac{z- w }{1-\bar{w}z}\right|$$ and
$P(0)=D$. We call $x\in M(H^\infty)$ a {\em trivial part} if $P(x)=\{x\}$.
Let
$$
G=\bigcup\{P(x):x \in M(H^\infty), \pee(x)\neq\{x\}\}.
$$
Then $G$ is an open subset of $M(H^\infty)$.

By Hoffman's work \cite{Hoffman}, $Z(b)\subset G$ for every
interpolating Blaschke product $b$ and for each $x$ in $G$, there exists
an interpolating Blaschke product $b$ such that $b(x)=0$. Also by
\cite{Hoffman}, for each $x\in G$ there exists a one-to-one and onto map $L_x:
D \to \pee(x)$ such that $f\circ L_x \in H^\infty$ for every $f\in H^\infty$.
The map $L_x$ is given as follows. Let $\{z_\alpha\}$ be a net
in $D$ with $z_\alpha \to x$ and let $L_{z_\alpha}(z) = (z +
z_\alpha)/(1+ \bar{z}_\alpha z)$. Then  $$(f\circ L_x)(z) =
\lim_\alpha (f\circ L_{z_\alpha})(z) \text{ for }f\in H^\infty \text{ and } z \in D.$$ 

A Blaschke product $b$ is of {\em type} $G$ if it is interpolating and $\{|b| <
1\} \subset G$. It is of {\em finite type} $G$ if it is of type $G$ and
for every $x\in Z(b)$ the set $Z(b)\cap P(x)$ is finite. A Blaschke
product is {\em locally thin} if for each $x\in Z(b)$ there
is an interpolating Blaschke product $q$ such that
$$
\lim_\alpha (1-|z_{n_\alpha}|^2) |\/q'(z_{n_\alpha})\/|=1
$$
whenever $\{z_\alpha\}$ is a subnet of the root sequence 
$\{z_n\}$ of $q$ that converges to $x$. Note $q$ may be different from
$b$. In fact, by \cite{Go-Li-Mo}, if $b=q$ for every $x\in
z(b)$, then $b$ is a thin Blaschke product.
Blaschke products of type $G$, finite type $G$ and locally thin Blaschke products are very important in the studies of Douglas algebras (see for example, \cite{Go-Li-Mo}, \cite{Gu}, \cite{Gu-Iz-I} and \cite{Gu-Iz-II}).

Let $\Omega$ be the family of all interpolating Blaschke products of finite type $G$ and $A$ be the closed subalgebra of $H^\infty$ generated by
$\Omega$. Let
$B=[A,\bar{A}]$ be the smallest (closed) $C^*$-subalgebra of
$L^\infty$ containing $A$. That is, $B$ is generated by the ratio of
interpolating Blaschke products of type $G$. Then $E=[H^\infty,
\bar{A}]$ will be the Douglas algebra generated by $H^\infty$ and
the complex conjugate of elements of $\Omega$.

Our main results are that every inner function $u$ in $E$ is a
finite product of interpolating Blaschke products of type $G$, 
from which we are able to identify the closure of the interpolatingBlaschke products of type $G$
among the Blaschke products. 
 As a consequence, we get $B=C_E$, where $C_E$ denotes the $C^*$-subalgebra of $L^\infty$ generated by the
invertible inner functions in $E$ and their complex conjugates. 
Another consequence is that the finite products of interpolating Blaschke
products of type $G$ are the only inner functions that are in $B\cap
H^\infty$. Hence by Theorem~4.1 of \cite{Ch-Ma}, $A=B\cap
H^\infty$ and finite convex combinations of finite products of
interpolating Blaschke products of type $G$ are dense in the closed
unit ball of $A$. For Blaschke product of finite type $G$, we obtain
similar results. In obtaining these results, we follow the approach in \cite{He}, but our proofs
rely heavily on the results about type $G$ and finite type $G$ developed in \cite{Gu-Iz-I} and \cite{Gu-Iz-II}.

Both authors would like to thank MSRI for support and hospitality.

\section{Results For Type $G$} 

 We begin with a few useful propositions concerning basic properties of interpolating Blaschke product of type $G$.

\begin{proposition}
If $B$ is of type $G$ and $b$ is a subproduct of $B$, then $b$ is of type $G$. If $b_1, b_2$ are of type $G$ and $b_1 b_2$ is an interpolating Blaschke product, then $b_1 b_2$ is of type $G$. 
If $b$ is of type $G$ and $b_{\lambda}=(b-\lambda)/(1-\bar{\lambda}b)$ is an interpolating Blaschke product, then $b_\lambda$ is of type $G$. The statements also hold if type $G$ is replaced by finite type $G$.
\end{proposition}
\begin{proof}
For type $G$, the first statement follows from $\{|b|<1\}\subset \{|B|<1\}\subset G$. The second statement follows from $\{ |b_1 b_2|<1\}=\{|b_1|<1\}\cup \{|b_2|<1\}\subset G$. The third statement follows from $\{|b_{\lambda}|<1\}=\{|b|<1\}\subset G$. 

For finite type $G$, the first statement follows from $Z(b)\cap P(m)\subset Z(B)\cap P(m)$. The second statement follows from $Z(b_1 b_2)\cap P(m)=(Z(b_1)\cap P(m))\cup (Z(b_2)\cap P(m))$. The third statement follows from Theorem 3.2 (iii) of \cite{Gu-Iz-II} because $H^\infty [\bar{b_{\lambda}}]=H^\infty [\bar{b}]$ by considering their maximal ideal spaces.
\end{proof}

We remark that not all of the statements in Proposition~1 are true for the family of thin Blaschke products.

\begin{proposition}
Suppose $b$ is an interpolating 
Blaschke product of type $G$ with roots $\{z_n\}$ in $D$. Let $q$ be
an interpolating Blaschke product with roots $\{w_n\}$ in $D$ such
that $\blip (w_n, z_n) \leq \curly  $ for all $n$ and for some
$\curly < 1$. Then $q$ is of type $G$.
\end{proposition}
\begin{proof}
Suppose $0< \lambda < 1$ and $z\in D$ such that
$|q(z)| < \lambda$. Then, by Lemma~1.4 and Corollary~1.3 on page 4 of
\cite{Garnett},
$$
|b(z)| = \prod^\infty_{n=1} \blip (z,z_n) \leq \prod^\infty_{n=1} 
\Bigl(\frac{\blip (z, w_n) + \curly }{1+ \curly \blip(z,w_n)}\Bigr)
       \leq \frac{|q(z)| + \curly }{1+ \curly |q(z)|}      
         < \frac{\lambda + \curly}{1 + \curly \lambda}
         = \lambda' < 1.
$$
So we have
$$
\{|q|<1\}=\bigcup_{0<\lambda<1} \{|q|<\lambda\}\subset \bigcup_{0<\lambda' <1} 
\{|b|<\lambda'\}=\{|b|<1\}\subset G.
$$
\end{proof}

\begin{proposition} Let $\cal{F}=\{ x\in M(H^\infty) : x \text{ is in 
the closure of some interpolating }$ $\text{sequences in } D \text{ whose Blaschke products are of type }
$G$ \}$. Then $\cal{F}$ is the union of a family of nontrvial Gleason parts.
\end{proposition}

\begin{proof}
By the definition of $\cal{F}$, every point in $\cal{F}$ belongs to a
nontrivial Gleason part. So let $m_0 \in \cal{F}$ and $m \in P(m_0)$.
Then $m_0 \in \overline{\{z_n\}}$ for some interpolating sequence
$\{z_n\}$ whose Blaschke product $b(z)$ is of type $G$. So there
is a subnet $\{z_{\alpha} \}$ of $\{z_n\}$ converging to $m_0$. Since $m\in
P(m_0)$, there is $\xi \in D$ such that
$$
\lim_\alpha L_{z_{\alpha}}(\xi) = L_{m_0}(\xi) = m.
$$
Let 
$$
\xi_n = L_{z_n}(\xi) = \frac{\xi + z_n}{1+ \bar z_n \xi}.
$$
Then $\blip (\xi_{n}, z_n) = |\xi| < 1$ for all $n$. By Corollary~1.6
on page 407 of \cite{Garnett}, there is a factorization $b=
b_1 b_2\cdots b_k$ with 
$$\delta(b_\jay) >
\frac{2|\xi|}{1 + |\xi|^2} \text{ for } \jay =1,2,\ldots, \kappa.$$ By
Lemma~5.3 on page 310 of \cite{Garnett},  each $Z(b_\jay)\cap D =
\{z_{\jay,n}\}$  is
interpolating. 

Since 
$$
\overline{\{z_n\}} = \bigcup^\kappa_{\jay =1}\overline{Z(b_\jay)}
$$
and the $Z(b_\jay)$'s have disjoint closures
\cite[p.~422]{Garnett}, it follows  that
$
m_0 \in \overline{Z(b_\jay)} = \overline{\{z_{\jay,n}\}}
$
for some $\jay$ and
the net $\{z_{\alpha}\}$ is eventually in $Z(b_\jay)$
because 
$$
M(H^\infty) \setminus \bigcup_{i\neq\jay}\overline{Z(b_i)}
$$ 
%
%[Silvio: check this formula] 
%
is an open neighborhood of $\overline{Z(b_\jay)}$ 
and $m_0\in \overline{Z(b_\jay)}$.

By Proposition~1, each $b_j(z)$ is of type $G$.
 By Proposition~2 (and the above argument), the Blaschke
product with roots $\{\xi_{\jay, n}\}$ is of type $G$. Finally,
$$
\lim_\alpha \xi_{\jay,\alpha} = \lim_\alpha 
L_{z_{\jay,\alpha}}(z) = L_{m_0}(\xi) = m
$$
and our assertion follows.
\end{proof}

\begin{proposition} An interpolating Blaschke product $b$ of type $G$
has modulus $1$ on those Gleason parts of $M(H^\infty)$ that do not
contain a point in $Z(b)$.
 Since $Z(b)$ equals the closure of $Z(b)\cap D$ and
$\cal{F}$ is the union of a family of Gleason parts, $b$ has in
particular modulus $1$ on $M(H^\infty)\setminus \cal{F}$.
\end{proposition}
\begin{proof} 
By Lemma~1.1 of
\cite{Gu-Iz-II} (or Theorem 1 of \cite{Gu-Iz-I})
we have $$\{|b|< 1\} = \bigcup_{m\in Z(b)}P(m).$$ Thus 
$$|b|=1 \text{ on } m(H^\infty)\setminus \bigcup_{m\in Z(b)}P(m).$$ The
second statement follows from the proof of Proposition 3.
\end{proof}

\begin{corollary} 
$M(E) = M(H^\infty)\setminus \cal{F}$.
\end{corollary}
\begin{proof}
By Theorem 1.3 on page 375 of \cite{Garnett}, 
$$
M(E)= \{m\in M(H^\infty): |b(m)|=1 \text{ for all } b\in \Omega \}.
$$
Now the results follows immediately from Proposition~4. 
\end{proof}

\begin{corollary}
$A$ is a proper subalgebra of $H^\infty$.
\end{corollary}

\begin{proof}
This follows because $M(H^\infty)\setminus (\cal{F}\cup M(L^\infty))$ is not empty.
\end{proof}

\begin{theorem} 
\label{thm:main}
Every invertible inner function in $E$ is a finite
product of interpolating Blaschke products of type $G$.
\end{theorem}
\begin{proof}
Let $u$ be an arbitrary invertible inner function in $E$, then $\bar{u}=u^{-1}$ in $E\subset L^\infty$. So $H^\infty[\bar{u}]\subset E$. By Corollary~5 above and Theorem 1.3 on page 375 of \cite{Garnett}
$$M(H^\infty)\setminus \cal{F}=M(E)\subset M(H^\infty[\bar{u}])=\{|u|=1\}.$$
Hence $\{|u|<1\}\subset \cal{F}\subset G$. This implies $Z(u)$ cannot contain any trivial parts. By Corollary~24 of \cite{McD-Su}, $u=b_1 b_2\cdots b_n$, where each
$b_j$ is an interpolating Blaschke product. Finally, for each $j$,
$$\{|b_j|<1\}\subset \bigcup_{k=1}^n \{|b_k|<1\}=\{|u|<1\}\subset G.$$
\end{proof}

\begin{corollary}
$B=C_E$, the $C^*$-subalgebra of $L^\infty$ generated by the inner functions
invertible in $E$ and their complex conjugates.
\end{corollary}

\begin{corollary}
Let $b$ be a finite product of interpolating Blaschke products of type
$G$. If $f\in H^\infty$ is such that
$\|f\|_\infty < 1$ and $\bar{f}b$
equals on $H^\infty$ a function $g$ almost everywhere on $T$, then the
function 
$$
b_f(z) = \frac{b(z) - f(z)}{1- g(z)},\quad\hbox{for } z\in D,
$$ 
is a finite product of interpolating Blaschke products of type $G$. 
\end{corollary}
\begin{proof}
Just observe that $b_f$ is an invertible inner function in $E$.
\end{proof}

In \cite{Ch-Ma}, Chang and Marshall showed that for an
arbitrary Douglas algebra $J$, the closed unit ball of $H^\infty \cap
C_J$ is the norm-closed convex hull of the Blaschke products in
$H^\infty\cap C_J$, where $C_J$ is the $C^*$-subalgebra of $L^\infty$ generated by
the invertible inner functions in $J$ and their complex conjugates. They
also showed that $J=H^\infty + C_J$ and that $D$ is dense in the maximal
ideal space of $H^\infty\cap C_J$. In our case $J=E$, $C_J=B$ and we have
the following corollary. (Note that an inner function in $B \cap H^\infty$ is
invertible in $E$.)

\begin{corollary}
\begin{enumerate} 
\item[(a)] $A=B\cap H^\infty$, and finite convex
combinations of finite products of interpolating Blaschke products of type $G$ are dense in the
closed unit ball of $A$.
\item[(b)] $E = H^\infty + B$.
\item[(c)] $D$ is dense in the maximal ideal space $M(A)$ of $A$.
\end{enumerate}
\end{corollary}

\section{Results For Finite Type $G$}

Next we will turn to the main results for interpolating
Blaschke product of finite type $G$ analogous to those established in section~2.
Let $\Omega_{\fin}$ be the family of all interpolating Blaschke products of finite type $G$ and $A_{\fin}$ be the
closed subalgebra of $H^\infty$ generated by 
$\Omega_{\fin}$. Let $E_{\fin} = [H^\infty,
\bar{A}_{\fin}]$, then $E_{\fin}$ is the Douglas algebra generated by $H^\infty$ and the
complex conjugate of elements of $\Omega_{\fin}$. 
We will 
show that Theorem~\ref{thm:main} holds if $E$ is replaced by $E_{\fin}$.

\begin{theorem}
\label{thm:was6}
Every invertible inner function in $E_{\fin}$ is a finite product of
interpolating Blaschke products of finite type $G$.
\end{theorem}

\begin{proof}
We will first show that the analog of Proposition~2 is true for
interpolating Blaschke products of finite type $G$.
Let $b$ be an interpolating Blaschke
product of finite type $G$ with zeros $\{z_n\}$ in $D$, and let $q$ be an
interpolating Blaschke product with zeros $\{w_n\}$ in $D$ such
that $\blip(z_n, w_n) \leq \curly $ for all $n$ and some 
$\curly < 1$. The proof of Proposition~2 shows that $\{|q|<1\}\subset
\{|b|<1\}$.  Since $b$ is of finite type $G$, by Theorem~2.1 of
\cite{Gu-Iz-II}, there is a subproduct $b_0$ of $b$ such
that $\{|b_0| < 1\} = \{|q| < 1\}$. 
 We will show that $q$ is of finite type $G$.

For $x\in Z(q)$, $|b_0(x)|<1$. By Lemma 1.1 of \cite{Gu-Iz-II}, there is an $x_0\in Z(b_0)$ such that $x\in P(x_0)$.
Suppose the set $Z(q)\cap P(x)=Z(q)\cap P(x_0)$ is infinite. Then, by
Theorem~3.1(i) of \cite{Gu-Iz-II}, there exist $y$ and
$y_0$ in $Z(q)$ such that $\supp \mu_y \subsetneq \supp \mu_{y_0}$.
Hence there are $m$ and $m_0$ in $Z(b_0)$ such that $y\in P(m)$ and
$y_0 \in P(m_0)$, but then $\supp \mu_m \subsetneq \supp
\mu_{m_0}$ (because by page 143 of \cite{Gamelin}, $\supp \mu_y = \supp \mu_m$ and $\supp \mu_{y_0} =
\supp  \mu_{m_0}$). Since $b_0$ is of finite type $G$, this contradicts Theorem~3.2(ii) of
\cite{Gu-Iz-II}. Thus $Z(q)\cap P(x_0)=Z(q)\cap P(x)$ must be finite.

Next we remark that the analogs of Propositions 3,4 and Corollaries 5,6  for finite type $G$ also hold by the same reasoning because of the analog of Proposition 2 for finite type $G$.
 
 Now let $u$ be an invertible inner function in $E_{\fin}\subset E$. By Theorem~7, $u=u_1 u_2\cdots u_m$, where each $u_i$ is of type $G$. Observe that
 if $\cal{F}_{\fin}$ is the analog of $\cal{F}$ for finite type $G$, then
 $$Z(u_i)\subset\{|u|<1\}\subset \cal{F}_{\fin}=\bigcup_{b\in \Omega_{\fin}} \{|b|<\frac{1}{2}\}.$$
Let $\delta_i=\inf\{\blip(w,z): w,z\in Z(u_i)\cap D, w\neq z\}>0$. Since $Z(u_i)$ is compact,
$$Z(u_i)\subset \bigcup_{j=0}^{n_i} \{|b_j|<\frac{1}{2}\},$$
for some $b_1, b_2, \ldots, b_{n_i}$ of finite type $G$. Let
$$S_{ij}=Z(u_i)\cap \{|b_j|<\frac{1}{2}\}\cap D,$$
then $$Z(u_i)\cap D=\bigcup_{j=1}^{n_i} S_{ij}.$$
By removing overlapping elements, we may assume the $S_{ij}$'s are disjoint.

Since $b_j$ is of type $G$, by Lemma 2.1 of \cite{Gu-Iz-II}, there is $\delta<1$ such that
$$S_{ij}\subset\{|b_j|<\frac{1}{2}\}\subset\{ z\in D:\blip(z,\{z_{j,n}\})\leq \delta\},$$
where $\{z_{j,n}\}$ is the root sequence of $b_j$ in $D$. For each disk $B(z_{j,n},\delta)=\{z\in D :\blip(z, z_{j,n})\leq \delta\}$, there is at most $k_i$ elements of $S_{ij}$ in $B(z_{j,n},\delta)$, where $k_i$ depends only on $\delta_i$.
So $S_{ij}$ is the union of at most $k_i$ sequences, each of which has at most one element in each $B(z_{j,n},\delta)$. By the analog of Proposition~2 for finite type $G$, proved above, the Blaschke product with root sequence $S_{ij}$ is a product of at most $k_i$ interpolating Blaschke product of finite type $G$. So $u_i$ is a 
finite product of at most $n_i k_i$ interpolating Blaschke product of finite type $G$. Therefore, $u$ is a finite product of interpolating Blaschke product of finite type $G$.
\end{proof}

In general, it is a difficult problem to determine the closure of an infinite set of interpolating Blaschke products among the family of all Blaschke products (see for example \cite{Li}). However, for Blaschke products of type $G$ and finite type $G$, their closures 
can be identified because of Proposition~1 and Theorems~7, 11.

\begin{theorem} Let $\cal{B}$ be the family of all Blaschke products with essential sup-norm. The closure of all interpolating Blaschke products of type $G$ in $\cal{B}$ is the set of all finite products of interpolating Blaschke products of type $G$. Also, the closure of all interpolating Blaschke products of finite type $G$ in $\cal{B}$ is the set of all finite products of interpolating Blaschke products of finite type $G$.
\end{theorem}

\begin{proof}
Suppose $B$ is in the closure of type $G$ products. Take $B$ of type $G$ such that $\|B-b\|_\infty=\|1-B\bar{b}\|_\infty<1$. It follows that $B\bar{b}$ is invertible in $E$ and so is $B=(B\bar{b})b$. By Theorem~7, $B$ is a finite product of interpolating Blaschke products of type $G$.

For the converse, it suffices to show for $B_1, B_2$ of type $G$ and $\varepsilon >0$, there is $B$ of type $G$ such that $\|B_1 B_2-B\|_\infty<\varepsilon$. Let the root sequences of $B_1$ and $B_2$ be $\{z_n\}$ and $\{w_n\}$, respectively. Since each of these sequences are separated, 
$$\delta_I=\inf \{\blip(z_m,z_n), \blip(w_m,w_n) : m\neq n\}>0.$$
By Hoffman's lemma (see Lemma 1.4 on pages 404-5 of \cite{Garnett}), there are 
$\delta_0, \varepsilon_0 <\delta_I/3$ such that
$$V_n\subset\{z\in D : \blip(z,z_n)<\varepsilon_0\},$$
$$W_n\subset\{z\in D : \blip(z,w_n)<\varepsilon_0\}$$
and $$\left\| z-\frac{z-(\delta_0/2)}{1-(\delta_0/2)z}\right\|_\infty<\varepsilon,$$
where $V_n$ and $W_n$ are the components of $\{z\in D: |B_1(z)|<\delta_0\}$ and $\{z\in D : |B_2(z)|<\delta_0\}$ containing $z_n$ and $w_n$, respectively. 
Since $3\varepsilon_0<\delta_I$, we have $$\blip(V_n,V_m)>\delta_I/3 \text{ and } \blip(W_n,W_m)>\delta_I/3 \text{ for} n\neq m.$$

Factor $B_2=B_3 B_4$ so that $w_n\in Z(B_3)$ if $\blip(w_n,Z(B_1))\geq \delta_0/4$ and $w_n\in Z(B_4)$ if $\blip(w_n,Z(B_1))<\delta_0/4.$ Let
$$B_5=\frac{B_4-(\delta_0/2)}{1-(\delta_0/2)B_4}$$ 
and $B=B_1 B_3 B_5$, then $\|B_1 B_2-B\|_\infty=\|B_4-B_5\|_\infty<\varepsilon$.

Next we will show $B$ is an interpolating Blaschke product. Since $\blip(Z(B_1),Z(B_3))\geq \delta_0/4$, $B_1 B_3$ is an interpolating Blaschke product. 
By Lemma 1.4 on pages 404-5 of \cite{Garnett}, $B_5$ is an interpolating Blaschke product. 

To see $\blip(Z(B_1 B_3),Z(B_5))>0$, let $w\in Z(B_5)$ and $z\in Z(B_1 B_3)$.  Then $w\in W_n$ for some $n$ and $w_n\in Z(B_4)$. If $m\neq n$, then
$$\frac{\delta_0}{2}=\blip(B_5(w),B_5(w_n))\leq \blip(w,w_n)<\varepsilon_0<\blip(w,w_m).$$
Since $w_n\in Z(B_4)$, we have $\blip(w_n,z_k)<\delta_0/4$ for some $z_k\in Z(B_1)$.

In the case $z\in Z(B_1)$ and $\blip(z,Z(B_4))<\varepsilon_0/4$, there is $w_m\in Z(B_4)$ such that $\blip(z,w_m)<\delta_0/4$ and so
$$\blip(w,z)\geq \blip(w,w_m)-\blip(w_m,z)\geq \frac{\delta_0}{2}-\frac{\delta_0}{4}=\frac{\delta_0}{4}.$$

In the case $z\in Z(B_1)$ and $\blip(z,Z(B_4))\geq \delta_0/4$, we have $z\neq z_k$, hence $\blip(z,z_k)\geq \delta_I$. So
$$\blip(w,z)\geq \blip(z,z_k)-\blip(z_k,w_n)-\blip(w_n,w)\geq\delta_I-\frac{\delta_0}{4}-\varepsilon_0\geq\frac{\delta_I}{2}.$$

In the case $z\in Z(B_3)$, we have
$$\blip(w,z)\geq \blip(z,w_n)-\blip(w_n,w)\geq \delta_I-\varepsilon_0\geq\frac{2\delta_I}{3}.$$
So $\blip(Z(B_1 B_3),Z(B_5))>0$. Therefore, $B=B_1 B_3 B_5$ is an interpolating Blaschke product. 

Finally since $B_1, B_2$ are of type $G$, by Proposition~1, $B$ is of type $G$. This completes the proof of the statement for the closure of type $G$. For the closure of finite type $G$, use Theorem~11 instead of Theorem~7 and repeat the above proof.
\end{proof}

\section{Questions}

Let $A_{\loc}$ be the closed subalgebra of $H^\infty$
generated by the locally thin Blaschke product and let $E_{\loc} =
[H^\infty; \bar{A}_{\loc}]$. 
\begin{enumerate}
\item[1.] Does Theorem~\ref{thm:main} hold if $E$ is replaced with $E_{\loc}$?
\item[2.] Does $E=E_{\loc}$ or does $E_{\loc} = E_{\fin}$, or neither?
\end{enumerate}
If we let $A^*$ be the closed subalgebra of $H^\infty$  generated by
the thin Blaschke products, and set $E^* = [H^\infty; \bar{A}^*]$,the main 
result of Hedenmalm 
\cite[Theorem~2.6]{He} asserts that Theorem~\ref{thm:main} holds for $E^*$. By
Proposition~1.1 and Example~3.1 of \cite{Gu-Iz-II}, there exists a Blaschke product $b$ of finite type $G$, which is not a finite product of thin Blaschke products.
It follows that $E^* \subsetneq E_{\fin}$ because otherwise $E_{\fin}=E^*$ would contain $H^\infty [\bar{b}]$ forcing $b$ to be a finite product of thin Blaschke products by \cite[Theorem~2.6]{He}. 
Also, by the proof of Theorem~2 of
\cite{Gu-Iz-I}, there exists a Blaschke product of type $G$, but not of finite type $G$. So we have $E_{\fin} \subsetneq E$. By Lemma~1 of \cite{Gu}, it can be shown
that $E_{\fin} \subseteq E_{\loc}$, but it is not clear whether $E_{\loc}
\subseteq E$ or $E_{\loc} \supseteq E$.

\end{document}